(This paper is in press in Duke Math. J. The substance
of the text is unaltered; inor changes and
corrections have been
made to the file to correspond to the version that will be published.)

\documentstyle[12pt]{article}
\input amssym.def
\input amssym.tex

\newcommand{\nc}{\newcommand}

\newcommand{\Hom}{\,{\rm Hom}\,}

\newcommand{\bla}{\phantom{bbbbb}}
\newcommand{\onebl}{\phantom{a} }
\newcommand{\eqdef}{\;\: {\stackrel{ {\rm def} }{=} } \;\:}

\newcommand{\half}{ {\frac{1}{2} } }

\newcommand{\beq}{\begin{equation}}
\newcommand{\eeq}{\end{equation}}
\newcommand{\beqst}{\begin{equation*}}
\newcommand{\eeqst}{\end{equation*}}
\newcommand{\barr}{\begin{array}}
\newcommand{\earr}{\end{array}}
\newcommand{\beqar}{\begin{eqnarray}}
\newcommand{\eeqar}{\end{eqnarray}}
\newtheorem{theorem}{Theorem}[section]

\newtheorem{prop}[theorem]{Proposition}
\newtheorem{definition}[theorem]{Definition}
\newtheorem{remit}[theorem]{Remark}

\nc{\lrar}{\longrightarrow}

\newcommand{\RR}{{\Bbb R }}

\newcommand{\ZZ}{{\Bbb Z }}

\newcommand{\cala}{{\mbox{$\cal A$}}}

\newcommand{\cali}{{\mbox{$\cal I$}}}

\newcommand{\calm}{{\mbox{$\cal M$}}}

\newcommand{\calp}{{\mbox{$\cal P$}}}

\newcommand{\calr}{{\mbox{$\cal R$}}}

\def\a{\alpha}
\def\b{\beta}
\def\g{\gamma}

\def\e{\epsilon}

\def\k{\kappa}
\def\l{\lambda}

\def\L{\Lambda}

\nc{\yyb}{{Y_\beta} }
\nc{\pr}{\partial}
\nc{\cbi}{C^{*,*} (G) }
\nc{\ctbi}{ {\tilde C}^{*,*}(G) }
\nc{\YY}{Y}
\nc{\mc}{\alpha}
\nc{\bmc}{\bar{\alpha} }
\nc{\liek}{  \mbox{\bf k}  }
\nc{\lieg}{  \mbox{\bf g}  }
\nc{\lieh}{  \mbox{\bf h}  }
\nc{\liehs}{ {\lieh}^* }
\nc{\lieks}{ {\liek}^* }
\nc{\liegs}{ {\lieg}^* }

\nc{\symf}{\omega}
\nc{\inpr}[1]{{ \langle {#1} \rangle} }
\nc{\epsr}{{ \epsilon_R} }
\nc{\zt}{{Z_t} }
\nc{\teta}{{ \tilde{\eta} } }
\nc{\mg}{{ \calm}}
\nc{\Proof}{ \noindent{\em Proof:} }

\nc{\diag}{\bigtriangleup}
\nc{\epm}{e}
\nc{\epmb}{{e_\beta}}

\nc{\sigfo}{\sigma}
\nc{\sigf}{\sigma}
\nc{\tw}{ {\tilde{ \omega}} }
\nc{\hmac}{K}
\nc{\hop}{I_\hmac}

\newcommand{\renorm}{{ \setcounter{equation}{0} }}

\nc{\FF}{ {\Bbb F} }
\nc{\tsig}{\tilde{\Sigma} }

\nc{\fund}{\Pi}
\nc{\free}{\FF}

\nc{\mb}{\calm_{\beta} }

\nc{\cpun}{ { \calp } }
\nc{\hatk }{ {K_c}}
\nc{\hk}{H^*_K}
\nc{\ttt}{t}
\nc{\bng}{\overline{NG} }
\nc{\ng}{NG}
\nc{\proj}{{\rm pr} }
\nc{\simp}{\bigtriangleup}

\nc{\phit} { {\tilde{\phi} } }
\nc{\xb}{X_\beta}

\nc{\ee}{{ \Bbb E} }
\nc{\pun}{E}
\nc{\bb}{{ \Bbb B} }
\nc{\psie} {{\Psi_E}}
\nc{\psib}{{\Psi_B} }
\nc{\bnk}{\overline{N \hatk } }
\nc{\nk}{{N\hatk} }
\nc{\uu}{\Bbb U}
\nc{\tfr} { {\tilde{f}_r} }
\nc{\tbrj}{ { {\tilde{b}_r}^j} }
\nc{\tar}{\tilde{a}_r}
\nc{\crk}{(c_r)_K }
\nc{\sio}{{\Sigma_0}}
\nc{\eex}{{ \Bbb E}_X }
\nc{\bbx}{{ \Bbb B}_X }
\nc{\ceex}{\check{ \Bbb E}_X }
\nc{\cbbx}{\check{ \Bbb B}_X }
\nc{\psiex}{ {\Psi_{X,E} }}
\nc{\psibx}{{\Psi_{X,B} } }
\nc{\cpsiex}{ {\check{\Psi}_{X,E} }}
\nc{\cpsibx}{{\check{\Psi}_{X,B} } }
\nc{\BPhi}{ \overline{\Phi} }
\nc{\BPhih}[1]{ \overline{\Phi^{\hmac}_{#1} } }
\nc{\Phih}{ \Phi^{\hmac}  }
\nc{\BPhihh}{ \overline{\Phi^{H} } }
\nc{\Phihh}{ \Phi^{H}  }

\begin{document}

\title{Group cohomology construction of the
 cohomology of  moduli spaces of flat connections
on 2-manifolds}

\author{Lisa C. Jeffrey \\
Mathematics Department \\
Princeton University \\
Princeton, NJ 08544, USA
\thanks{This material is based on work
supported by the National Science Foundation under Grant No.
DMS-9306029.} }

\date{April  27, 1994}

\maketitle

\begin{abstract}
We use group cohomology and the de Rham complex on simplicial manifolds
to give explicit differential forms representing generators of the
cohomology rings of moduli spaces of representations
of fundamental groups  of 2-manifolds. These generators are constructed using
the de Rham representatives for the cohomology of classifying spaces
$BK$ where $K$ is a compact Lie group; such representatives (universal
characteristic classes) were found by Bott and Shulman.  Thus our
representatives for the generators of the cohomology of moduli spaces
are given explicitly in
terms of the Maurer-Cartan form.
This work solves a problem posed by
Weinstein, who gave a corresponding construction (following Karshon
and Goldman) of the
symplectic forms on these moduli spaces.

We also give a corresponding construction of equivariant differential
forms on the extended moduli space $X$, which is a  finite dimensional
symplectic space equipped with a Hamiltonian action of $K$ for which
the symplectic reduced space is the moduli space of representations
of the 2-manifold fundamental group in $K$.

\end{abstract}

\renorm
\section{Introduction} \label{s1}

Let $K$ be a compact connected semisimple
Lie group, and let $\b$ be an
element in the center $Z(K)$.
 If $\Sigma$  is a closed 2-manifold
of genus $g \ge 2$ with fundamental group $\fund$,
and $\sio = \Sigma - D^2$, then   $\pi_1(\sio) = \free$ where
$\free $ is the free group on $2g$ generators.  There is an associated
moduli space of representations
 $\calm_\b = Y_\b/K$ where $\yyb = \{ \rho \in {\rm Hom} (\free, K):
\rho(R) = \b \} $ and $K$ acts on ${\rm Hom} (\free, K) $ by conjugation.
Here $R$ is the element of $\free $ corresponding to a loop winding
once around the boundary $\partial \Sigma_0$. Because $\beta $ is in the
center $Z(K)$, each point in the space $\calm_\beta$ gives rise to
a representation of the fundamental group $\fund$ of the closed
surface $\Sigma$ into the group $K_c = K/Z(K)$.

The space $\calm_\b$ has two alternative descriptions. Via the
holonomy map, $\calm_\b$ may be identified with the space of gauge
equivalence classes of flat connections on a principal
$K$ bundle over $\sio$, for which the  holonomy around the
boundary $\partial \sio$ is the
element $\beta$.  We obtain a second alternative description
once we fix a complex structure on $\Sigma$: then $\calm_\b$
becomes identified with a  space of semistable holomorphic
vector bundles (of prescribed rank and degree) over $\Sigma$.

Atiyah and Bott worked in this holomorphic setting, and
described the generators of the cohomology ring of
$\calm_\b$ in terms of a holomorphic vector bundle
$\Bbb U$ over $\calm_\b \times \Sigma$ (the {\em universal
bundle}), whose restriction to a point $m \in \calm_\b$ is the
holomorphic vector bundle over $\Sigma$ corresponding to $m$.
The K\"unneth decomposition of the Chern classes of $\Bbb U$
yields classes in $H^*(\calm_\b)$ which are the generators of the
cohomology ring. The purpose of the present paper is to give an
explicit description of these generators in a representation theoretic
setting, using group cohomology.

Our starting point is the paper [W] of Weinstein: the construction
we present below generalizes Weinstein's construction of the symplectic
form on moduli space, whose cohomology class is one of the generators
described above. Goldman [G] constructed the symplectic form using
group cohomology, but his proof that this form was closed used
its gauge theory description. Karshon [K] gave the first group
cohomology proof that the symplectic form was closed. In
[W], Weinstein  interpreted Karshon's construction in terms of the
realization of $H^*(BK)$ via the de Rham cohomology of simplicial
manifolds, which is due to Bott and Shulman [B1,Sh]. Here we extend
Weinstein's treatment to all the generators of
$H^*(\calm_\b)$;\footnote{In this paper,
all cohomology groups are taken with complex coefficients; one could
give an equivalent treatment using rational or real coefficients.}
this problem was posed in the final section of [W].

The key advance which facilitates our construction is a description
of the universal bundle and its base space $\Sigma\times \yyb$ as
simplicial manifolds; such a description  is possible since
$\Sigma$ is a realization of $B \fund$ and $\sio$ a realization of
$B \free$. The spaces $EK$ and $BK$ may also be described as
simplicial manifolds, and we exhibit the classifying map for
the universal bundle as a map of simplicial manifolds. The generators
of $H^*(\calm_\b)$ are given by pairing the characteristic classes
of the universal bundle with homology classes of $\Sigma$ via the slant
product; in the present setting the relevant homology
classes are classes in the group homology of $\fund$.

One feature which  we generalize from [W]  is the role of
equivariant cohomology. For any group $G$,
the total space of the classifying bundle $EG$ over the
classifying space $BG$ is equipped with an action of $G$ which
covers the action of $G$ on $BG$ associated to the
adjoint action of $G$ on itself.  We construct a universal
bundle $\ee$ over $\Sigma \times \yyb$, which is equipped with an
action of $K$ covering the action of $K$ on $\Sigma \times \yyb$.
The classifying map for $\ee$ is a $K$-equivariant map; our
de Rham representatives for the generators of the cohomology
ring of $\calm_\b$ (Theorem \ref{t7.1})
are obtained by using the classifying map to pull back
{\em equivariant} cohomology classes
from the classifying space. (The
appropriate classifying space is $B \hatk$ rather than $BK$, where
 $\hatk = K/Z(K)$.)  These equivariant characteristic classes appear
naturally in the {\em Cartan model} for equivariant cohomology via
an equivariant version of the Chern-Weil map [BGV]; we evaluate
invariant polynomials on an equivariant extension of the curvature form
of a distinguished  invariant connection on the classifying bundle
$E \hatk$
over $B \hatk$. Thus we have constructed the \lq grand unified theory'
conjectured in the final section of [W], in which extensions
of characteristic forms to equivariantly closed forms appear naturally.

Finally in Theorem
\ref{t8.1} we generalize our construction to give equivariantly closed forms on
the {\em extended moduli space} $\xb$, a noncompact Hamiltonian
$K$-space which contains
$\yyb$ as the zero locus of the moment map,
so that  $\calm_\b$ may be obtained from  $\xb$ by symplectic reduction.
The explicit description of the
de Rham representatives
for the generators
(Propositions \ref{p7.*2} and
\ref{p8.2})
 will be a key step in a future paper on intersection pairings in
$H^*(\calm_\b)$, where we shall require explicit representatives for the
equivariant cohomology classes on $\xb$, to estimate their
growth and show they satisfy appropriate
boundedness properties as the value of the moment map on $\xb$ tends
to infinity.

This paper is organized as follows. In Section \ref{s2} we give the
representation theory definition of the universal bundle, which is reworked
in the setting of simplicial manifolds in Section \ref{s6}. Section \ref{s3}
contains background material on simplicial manifolds and the de Rham
complex on simplicial manifolds;
Section \ref{s4} explains the simplicial realization of the classifying
space $BG$ and the Bott-Shulman construction of universal characteristic
classes, while Section  \ref{s5} gives the generalization
of the Bott-Shulman construction to equivariant cohomology.
Sections \ref{s7} and \ref{s8} contain our final results:
Section \ref{s7} exhibits the generators of the cohomology ring
of $\calm_\b$ by pulling back  equivariant characteristic
classes from $H^*(B \hatk)$ via the classifying map, while
Section \ref{s8} gives the generalization to equivariant cohomology
classes on the extended moduli space $\xb$ which restrict on
$\yyb$ to the classes from Section \ref{s7}.

\noindent{\em Acknowledgement:} We thank Alan Weinstein for
helpful conversations, and Ezra Getzler for an explanation of the
results of [Ge] and their relevance to this paper.

\renorm
\section{Definition of the universal bundle}  \label{s2}

Let $K$ be a compact connected semisimple
Lie group, and let
$\hatk = K/Z(K)$.
 We start with the following
construction of a   family of principal bundles over a
closed 2-manifold $\Sigma$ of genus $g \ge 2$ (the
universal bundle)  which generalizes
a construction given in [T]. Let $\tsig \cong D^2$ be the
universal cover of $\Sigma$, and let $\fund = \free/\calr$ be the
fundamental group of $\Sigma$; here, $\free = \free_{2g}$ is the free group on
$2g$ generators $x_1, \dots, x_{2g}$, and $\calr$ is the normal subgroup
generated  by the relator $R = \prod_{j = 1}^g x_{2j-1} x_{2j} x_{2j-1}^{-1}
x_{2j}^{-1}$. If $\b$ is a prescribed element of the center $Z(K)$, we
let
\beq \label{1.1} \yyb = \{ \rho \in \Hom(\FF, K) \, : \,
\rho(R) = \beta \}. \eeq
Define $\calm_\b = \yyb/K$; for some groups $K$ and some choices of $\b$
(for instance $K = SU(N)$ and $\b$ a generator of $Z(K)$) the space
$\calm_\b$ is a smooth manifold.
Instead of working on $\calm_\b$ we shall equivalently work with
$K$-equivariant objects on $\yyb$.

We form the  universal bundle $\pi: \pun \to $  $B = \Sigma \times \yyb$
as follows, with fibre $\hatk = K/Z(K)$.
The universal bundle may equivalently be described as a family
$\pun$ of  flat bundles over $\sio = \Sigma - D^2 $, such that for a point
$\rho \in Y_\b$  (for which  the corresponding point $[\rho] \in $
$\calm_\beta$ specifies  an equivalence class of flat
bundles), the bundle $\pun|_{\sio \times \rho} $ is the flat
bundle parametrized by $\rho$.
We define the bundle $\pun$ by
\beq \label{2.02} \pun = ( \tsig \times \yyb \times \hatk)/\fund  \eeq
where $p \in \fund$ $= \pi_1(\Sigma)$
 acts\footnote{We write the
action of $K$ on $\hatk$ (by composing the quotient
map $K \to \hatk$ with the action of $\hatk$ on itself by left
multiplication) as
$t \in K: \onebl k \in \hatk \mapsto tk \in \hatk.$}
on
$ \tsig \times \yyb \times \hatk $ as follows:
\beq \label{1.2} p: \; (\sigma, \rho , k ) \to (\sigma p, \rho,
\rho(\tilde{p} )^{-1} k ). \eeq
Here, $\tilde{p} $ is an element of $\free$ which descends to $p \in \fund$
under the quotient map. Notice that since $\rho(R) = \beta \in Z(K)$,
the element $\rho(\tilde{p} ) $ is well  defined in $\hatk $ $= K/Z(K)$.
The map $\pi: \pun \to B = \Sigma \times \yyb$ is given by
$\pi (\sigma, \rho, k) = ([\sigma], \rho)$.

Now the bundle $\pun$ is equipped with an action of $K$, which is given by
an action on $\tsig \times \yyb \times \hatk$ as follows:
\beq \label{1.3} \ttt \in K: (\sigma, \rho, k)
\mapsto (\sigma, \ttt \rho \ttt^{-1}, \ttt k). \eeq
It is easy to verify that this action is compatible with the action
of $\fund$ on $\tsig \times \yyb \times \hatk$, and so it descends to a well
defined action on the bundle $\pun$.
The $K$-equivariant principal bundle $\pun$ over $\Sigma \times \yyb$
may  be thought of as a principal bundle over $\Sigma \times \calm_\b$.

\renorm
\section{Simplicial manifolds} \label{s3}
Let $\simp^n$ $ = \{ (t_0, \dots, t_n) \in
[0,1]^{n+1}: \sum_i t_i = 1 \} $ denote the standard $n$-simplex.
A {\em simplicial manifold} $X =
\{n \mapsto X(n) \} $ is a
contravariant functor from the category of simplices to the category
of $C^\infty$ manifolds. In particular, to every
nonnegative integer $n$ there is associated a manifold
 $X(n)$ (corresponding to $\simp^n$), and there are a
collection of maps
$\e_i: X(n) \to X(n-1) $ (the {\em face  maps}) for $0 \le i \le n$
which are functorially associated to the
inclusion maps  $\e^i: \simp^{n-1} \to \simp^n$ of the $i$'th face.
(There are also {\em degeneracy maps}  $\delta^*: X(n-1) \to X(n)$
associated to simplicial
maps $\delta: \simp^n \to \simp^{n-1}$, but these will
not directly concern us.)
See for instance  [Du,Ge,Se1,Se2] for
a more extensive discussion of simplicial manifolds, and Chapter  1 of [M]
for background on simplicial objects.

There are two versions of the de Rham complex on a simplicial
manifold $X$, which appear in the work of Dupont [Du] (Section 2).
One version corresponds to the {\em fat realization}
$\parallel X \parallel $ of $X$, which is defined
as the quotient space of  $\coprod_n \simp^n \times X(n)$
by the equivalence relation
\beq  (\e^i \times {\rm id}) (t,x) \sim
({\rm id} \times \e_i) (t,x), \bla (t, x) \in \simp^{n-1} \times X(n). \eeq
Following Dupont [Du], we then make the following definition
of a $k$-form on $X$:
\begin{definition} \label{d3.1}
A $k$-form on $X$ is a collection of $k$-forms
$\phi^{(n)} \in \Omega^k(\simp^n \times X(n) ) $ satisfying
\beq \label{3.2} ({\rm id} \times \e_i)^* \phi^{(n-1)}
= (\e^i \times {\rm id})^* \phi^{(n)} \eeq
on $\simp^{n-1} \times X(n)$ for all $0 \le i \le n$ and
all $n$ $\ge 1$. \end{definition}
\begin{definition} The bicomplex $(A^{*, *} (X), d_\simp, d_X) $
is the complex $\bigoplus_{k,l} A^{k,l} (X) $, where
$A^{k,l} (X)$ is the vector space spanned by
 forms $\phi =  \{ \phi^{(n)} \} $ on $X$
(in the sense of Definition \ref{d3.1}) for which $\phi^{(n)} $
is a linear combination of forms $ \omega_k \wedge \eta_l$
where $\omega_k \in \Omega^k(\simp^n) $ and
$\eta_l \in \Omega^l (X(n) ) $.
This bicomplex  is equipped with the differentials
$d_\simp$ and $d_X$ which are the
exterior differentials on $\simp^n$ and $X(n)$ respectively.
\end{definition}

There is a second bicomplex $\cala^{*,*} (X)$
which represents the de Rham
complex on the simplicial manifold $X$:

\begin{definition} The bicomplex $(\cala^{*,*} (X), \delta, d)$
is defined by
\beq \cala^{k,l} (X) = \Omega^l (X(k) ); \eeq
we define the differential  $\delta$  by
\beq \label{3.01}
\delta =  \sum_i (-1)^i \e_i^*,  \eeq
while  $d= d_X$ is the exterior differential on $X(k)$ for all $k$.
\end{definition}

Lemma 2.3 of [Du] gives a chain homotopy equivalence  between the
double complexes $A^{*,*} (X)$ and $\cala^{*,*} (X)$. In one
direction the
map $\cali: A^{k,l}(X) \to \cala^{k,l} (X) $ giving rise to
the chain homotopy is defined by  integration over the simplex $\simp^n$:
$$\cali: \phi^{(n)} \in \Omega^*(\simp^n \times X(n))
 \mapsto \int_{\simp^n}  \phi^{(n)}. $$

Let $G$ be a real Lie group with Lie algebra $\lieg$. Then there is a standard
realization  (see for example [Se1]) of the classifying space $BG $ and the
universal bundle $EG \to BG$ as  simplicial  manifolds, which is given as
follows. (See also [Du], Section 3.)  We define simplicial
manifolds $\bng$ and $\ng$ by
\beq \label{2.1} \bng(n) = G^{n+1} \eeq
\beq \label{2.2} \ng(n) = G^n \eeq
Here, the face maps $\bar{\e_i}: \bng(n) \to \bng(n-1) $
are as follows:
\beq \label{3.6a}  \bar{\e_i} (g_0, \dots, g_n)  = (g_0, \dots, \hat{g_i} ,
\dots, g_n)
\bla (i = 0, \dots, n) \eeq
while the face maps $\e_i: \ng(n) \to \ng(n-1) $ are given by
\beq \label{3.6b}
\begin{array}{lcr} \e_0 (h_1, \dots, h_n) &=& (h_2, \dots, h_n)  \\
\e_i (h_1, \dots, h_n) &=& (h_1, h_2, \dots, h_i h_{i + 1} , \dots,
h_n)  \bla (i = 1, \dots, n-1)  \\
 \e_n(h_1, \dots, h_n) &=& (h_1, \dots, h_{n - 1} ). \end{array}
\eeq
There is then
a map of simplicial manifolds $q: \bng \to \ng$ given by
a collection of maps $q_n: \bng(n) \to \ng(n)$,
defined by
\beq \label{2.3} q_n (g_0, \dots, g_n) = (g_0 g_1^{-1}, \dots,
g_{n-1} g_n^{-1}). \eeq
The simplicial $G$-bundle $\bng \stackrel{q}{\to} \ng$ is the
simplicial realization
 of the $G$-bundle $EG \to BG$. The action of
$G$ on the total space $\bng(n) $
of the principal $G$-bundle
$\bng(n) \to \ng(n)$  is given by
the action of $G$ on $\bng(n) = G^{n + 1} $ by right
multiplication.
In the case $X = \ng$, the
bicomplex $(\cala^{*,*} (X), \delta, d)$ is the bicomplex referred to
in
[W] as  the {\em de Rham-bar bicomplex}.

Now suppose the group $G$ is acted on by another group $H$, in other words
that there is a homomorphism from $H$ to ${\rm Aut} (G)$. Then the simplicial
manifold  $NG$ is also acted on by $H$. We shall consider the case when
$H=G$ acts on $G$ by the adjoint action. In this case the action of
$G$ on itself by   {\em left} multiplication gives rise to an action
of $G$ on $\bng$ which covers the adjoint action of $G$
 on $\ng$.\footnote{This action is not
to be confused with the action of $G$ on
$\bng$ by {\em right} multiplication, which gives $\bng $ its structure
of principal $G$-bundle over $\ng$.}
Explicitly, the action of $G$ on $\bng(n) $ is given by
$$g \in G: (g_0, \dots, g_n) \mapsto (g g_0, \dots, g g_n), $$
while the action on $\ng(n) $ is given by
$$ g \in G : (\g_1, \dots, \g_n) \mapsto (g \g_1 g^{-1},
\dots, g \g_n g^{-1} ) . $$
 Thus for this
particular action (the adjoint action of $H = G$ on $G$),
we have exhibited the structure of an $H$-equivariant bundle on the
simplicial bundle $\bng \to \ng$. The
existence of this $H$-equivariant bundle structure  will be crucial for what
follows: it does not hold for general actions of $H$ on $G$.

\renorm
\section{Characteristic classes} \label{s4}
In this section,
let $G$ be a compact connected  Lie group with Lie algebra
$\lieg$.  In terms of the construction
of $\bng$ and $\ng$ as simplicial manifolds, there is a standard
construction of generators for $H^*(BG)$ (in other words, of characteristic
classes) due to Bott and Shulman [B1,Sh].
According to these references, an invariant polynomial
$Q \in S(\liegs)^G$ of degree $r$ gives
rise to differential
forms $\BPhi_n(Q) \in \Omega^{2r-n} (G^{n+1} ) $ which in fact
pull back from differential forms $\Phi_n(Q) \in
\Omega^{2r-n}(G^n)$ under the projection maps
$q_n$ defined in (\ref{2.3}).
 These forms are compatible with the face maps of $\bng$
and $\ng$, and fit together to form a closed form of degree $2r$
(in the sense of Definition \ref{d3.1}) on
$NG$ (see [Sh] Proposition 10 or [B1] Section 1).

These differential forms are constructed as follows (see  [Du] Lemma 3.8,
[Ge], and [Sh] Section II). We let $\theta$ $\in \Omega^1 (G) \otimes
\lieg$ be the (left invariant) Maurer-Cartan form on $G$.
Also, let $\bar{\theta}  \in \Omega^1(G) \otimes \lieg$
be the corresponding right invariant form: in other
words $\theta(g \eta) = \eta= \bar{\theta}(\eta g)$,
where $g \in G$ and $\eta\in \lieg$. Defining
$\proj_i: G^{n+1} \to G$ ($i = 0, \dots, n$) as the projection on
the $i$'th copy of $G$, we denote
$\proj_i^* \theta$ by $\theta_i$.

Let $t = (t_0, \dots, t_n) \in \simp^n$, where
$\simp^n = \{  (t_0, \dots, t_n) \in \RR^{n+1}: \onebl $
$ \sum_{i = 0 }^n t_i = 1 \} $ is the
standard $n$-simplex.
We define
 $\theta $ $ \in \Omega^1 (\simp^n \times G^{n+1} ) \otimes \lieg$ by
\beq \label{3.1}
\theta: t \in \simp^n \mapsto
 \theta(t) = \sum_{i = 0 }^n t_i \theta_i \in
\Omega^1(\simp^n \times G^{n+1} ) \otimes \lieg  \eeq
This is a connection on the principal
$G$-bundle $\simp^n \times \bng(n) \to
\simp^n \times \ng(n)$. Define
$$F_{\theta(t) } = d(\theta(t) ) + \half [\theta(t), \theta(t) ] $$
to be the curvature of the connection $\theta(t)$.
We then define
\beq \label{3.101} \BPhi_n(Q) = \int_{\simp^n}
Q(F_{\theta(t) } ) \in \Omega^{2r-n} (G^{n+1} ).  \eeq
In other words, $Q(F_{\theta(t)} ) $ represents an
element in $A^{*, *} (\bng) $ and $\BPhi_n(Q) $
represents the corresponding element in $\cala^{*, *}(\bng)$.
The $\BPhi_n(Q)$ give a differential form $\BPhi(Q)$
on the (simplicial) total
space $\bng$, which pulls back from a differential form  (in the
sense of Definition \ref{d3.1})
$\Phi(Q)  = \{ \Phi_n(Q) , n = 1, \dots, r\}$ on  $\ng$. One may take
$\Phi_n(Q) = \sigma_n^* \BPhi_n(Q)$, where
$\sigma_n$ are sections
 of the bundles $q_n: \bng(n) \to \ng(n)$
 (\ref{2.3}) given by
\beq \label{4.2a}
 \sigma_n:
(g_1, \dots, g_n) \mapsto (g_1 g_2 \cdots g_n, g_2 \cdots g_n,
 \dots, g_n, 1 ).
\eeq
Each term in $F_{\theta(t)} $ contains at most one power of $dt_i$, so
for $n > r $, $Q(F_{\theta(t)}) $ restricts to
$0$ on $\simp^n$ and one finds $\Phi_n(Q) = 0 $
([Sh], Proposition 10).
The map $Q \mapsto \Phi(Q) $ is called the {\em Bott-Shulman map}.
By construction $\Phi(Q)$ is closed under the total
differential $D$ which is given
(on elements of form  degree $p$ in $X = \ng$) by
 $D = \delta + (-1)^p  d$.\footnote{We use the following
(standard) sign convention for
total differentials of double complexes $C^{*,*} $ with differentials
$d_1: C^{p,q} \to C^{p+1,q} $ and
$d_2: C^{p,q} \to C^{p,q+1} $. If $d_1 d_2 = - d_2 d_1$  then the total
differential is $D = d_1 + d_2$. If $d_1 d_2 = d_2 d_1$
-- as
is the case here -- then the total
differential $D$ on $C^{p,q} $ is $d_1 + (-1)^p d_2$.}
In other words,
\beq\label{4.m1}
(\delta \pm d) (\Phi_1 (Q) + \Phi_2(Q) + \dots + \Phi_r(Q) ) = 0 \eeq
or
\beq \label{4.m2} d \Phi_1(Q	)  = 0 , \eeq
\beq \label{4.m3} \delta \Phi_{j - 1} (Q) = (-1)^{j+1}
 d \Phi_j(Q	)  \in
\Omega^{2r-j+1} (G^j)  \onebl (1 < j \le r),  \eeq
\beq \delta \Phi_r(Q) = 0. \eeq
The Bott-Shulman map
gives a collection of differential forms on
the simplicial realization $\ng$ of $BG$, which are representatives
in de Rham cohomology for the elements of $H^*(BG)$.

\noindent{\em Example.} Let $G$ be a compact connected Lie group
with Lie algebra $\lieg$, and let
\beq \label{4.001}Q_2  = \langle \cdot, \cdot \rangle \eeq
be the quadratic form on $\lieg$ determined by an
invariant inner product $\langle \cdot, \cdot \rangle$.
Then we find that
\beq \label{4.003} \Phi_1 (Q_2) = -  \frac{1}{6} \langle \theta,
[\theta, \theta] \rangle  \eqdef -  \lambda  \in \Omega^3(G), \eeq
\beq \label{4.004} \Phi_2(Q_2) = \langle \theta_1, \bar{\theta}_2 \rangle
\eqdef \Omega \in \Omega^2(G^2). \eeq
Here, $\theta_1  = \proj_1^* \theta$ is the
(left invariant) Maurer-Cartan form on the first copy
of $G$ and $\bar{\theta}_2 = \proj_2^* \bar{\theta} $ is the
corresponding right invariant form on the second copy of $G$.

\renorm
\section{Equivariant characteristic classes} \label{s5}

\noindent{\em  5.1 The equivariant simplicial de Rham complex}

Let us now consider the equivariant analogue of the construction
given in the previous section.
For a general manifold $M$ acted on by a compact connected  Lie
group $H$, there is an equivariant analogue $\Omega^*_H(M)$ of the
de Rham  complex  (the {\em Cartan model} -- see  [BGV], Section 7.1)
 whose cohomology is the equivariant cohomology
$H^*_H(M)$. The complex  $\Omega^*_H(M)$ is defined
by
$$ \Omega^*_H(M) = \Bigl ( \Omega^*(M) \otimes S(\lieh^*) \Bigr )^H. $$
The grading on this complex is given as follows. If
$\alpha: \phi \in \lieh  \mapsto \alpha(\phi) \in \Omega^*(M)$ is
an element of $\Omega^*_H(M)$ which is of degree $r$ as a polynomial
in $\phi$ and of degree $p$ as a differential form on $M$, then it is
assigned grading $p + 2r$. The differential on
$\Omega^*_H(M) $ is given by
\beq \label{5.01}
 (d_H \alpha) (\phi) = (d  - \iota_\phit ) (\alpha( \phi) ) \eeq
where $\phit$ is the vector field on $M$ generated by $\phi$
and $\iota$ is the interior product.

\noindent{\em Remark.} The operators $d$ and $\iota_\phit$ on
$\Omega^*_H(M)$ in fact anticommute (since  $  d \iota_\phit +
\iota_\phit d = L_{\phit} $, and the Lie derivative
$L_{\phit} $ vanishes on $\Omega^*_H(M)$.) Thus
$\Omega^*_H(M)$ may itself be given a bicomplex structure (see
e.g. [H], Section 1); we shall not need this structure here.

If a compact connected Lie group $H$ acts on
a  {\em simplicial}
 manifold $X$, we may form two double  complexes  which compute
the equivariant cohomology of $X$, by
analogy with the double complexes which appeared in Section \ref{s3}.
One such complex $(A^{*,*}_H(X), d_\simp, d_H)$ is defined as
\beq A^{*,*}_H(X) = (A^{*,*} (X) \otimes S(\liehs))^H,  \eeq
where  the differentials are the exterior differential $d_\simp$
on $\simp^n$ and the Cartan model equivariant exterior differential
$d_H$ (\ref{5.01}) on $X(n)$.
The other double complex $(\cala^{*,*}_H(X), \delta, d_H) $
is defined by
\beq \cala^{n,q}_H(X) = \Omega^q_H(X(n) )  \eeq where
the differential $\delta$ is still given by (\ref{3.01}) and the
differential $d_H$ is given by (\ref{5.01}). As in Section \ref{s3},
a map from $A^{*,*}_H (X) $ to $\cala^{*,*}_H (X)$ is given by
integration over $\simp^n$.

The cohomologies of the total complexes of $A^{*,*}_H(X) $
and $\cala^{*,*}_H(X) $  are isomorphic, since Dupont's construction
([Du], Theorem 2.3) of a chain  homotopy  equivalence between the
two complexes extends to the equivariant case. Applying minor modifications
of the proofs
 in [Ge]
(combining Theorem 2.2.3 with the evident extension of
Corollary 1.2.3 to simplicial manifolds), one finds that the
cohomology of either of
these total  complexes  indeed equals the equivariant
cohomology of the simplicial manifold $X$.

\noindent{\em 5.2 Equivariant characteristic classes}

Now let $X$ be a manifold acted on by a  compact connected
 Lie group $H$,
and let $P$ be a principal $G$-bundle over $X$
(for some compact connected Lie group $G$), such that $H$ acts
on the total space  of $P$ compatibly with its action on $X$. Let
$\theta \in \Omega^1(P) \otimes \lieg $ be a connection on $P$,
invariant under the action of $H$. We then define the {\em moment}
$\mu \in (C^\infty (P) \otimes \lieg)^H \otimes \lieh^* $ (see Section 7.1
of [BGV]) by
\beq \label{4.002} \phi \in \lieh \mapsto
\mu(\theta)  (\phi) = - \iota_\phit \theta \in \lieg, \eeq
where $\phit$ is the vector field on $P$ generated by $\phi \in \lieh$.
Given an invariant polynomial $Q \in S(\liegs)^G$, we may form the
corresponding equivariant characteristic class of $P$, which
has a  representative
in the Cartan model   given by ([BGV]
Theorem 7.7)
\beq  Q (F_\theta + \mu) \in
(\Omega^*(P) \otimes S(\lieh^*) )^H. \eeq
This differential form represents a class in
$H^*_H(P) $ which is the pullback of a class in
$H^*_H(X)$, the equivariant characteristic class
corresponding to the invariant polynomial $Q$.

\noindent{\em 5.3 Equivariant version of the  Bott-Shulman map}

We shall now specialize the construction of Section 5.1 to the case
when the simplicial manifold $X$ is $\bng$ or $\ng$, and the
group $H$
acts on $\bng$ and $\ng$ compatibly with the projection map from
$\bng$ to $\ng$. We require also that the connection
$\theta(t)$ (\ref{3.1})
on $\simp^n \times \bng(n)$ be preserved by the
action of $H$. This is true, for instance, when $G$ is a quotient of $H$
and $H$ acts via the
left action of $G$ on $\bng$ and
the adjoint action of $G$ on $\ng$.\footnote{We shall apply this when
$H = K$ is compact semisimple and $G = \hatk  = K/Z(K)$.}
The connection
$\theta(t) $ (\ref{3.1}) on
$\simp^n \times \bng(n) $ is then invariant under this action of $G$.
Hence, according to Section 5.2, an invariant polynomial
$Q$ on $\lieg$ of degree $r$ gives rise to an element
$$\BPhihh_\simp (Q)=
Q\Bigl (F_{\theta(t) } + \mu(\theta(t) ) \Bigr ) \in
A^{*,*}_H(\bng), $$
which pulls back from
$\Phihh_\simp (Q) \in A^{*,*}_H(\ng).$
Integration over the simplex $\simp^n$ gives an element
$$ {\BPhihh}(Q) =
\int_{\simp^n} Q\Bigl (F_{\theta(t) } + \mu(\theta(t) ) \Bigr )
 \in \cala^{*,*}_H(\bng), $$
which  pulls back from $\Phihh(Q) \in \cala^{*,*}_H(\ng) $
and
represents the equivariant characteristic class associated to
$Q$ in the Cartan model for  $\Omega^*_H(BG)$.
We denote the component on $\ng(n)$ by $\Phihh_n(Q)  \in
\Omega^{2r-n}_H(\ng(n) ) $.
In the same way as in (\ref{4.m1}-\ref{4.m3}), we have
\beq \label{5.m1}
(\delta \pm d_H) (\Phihh_1 (Q) + \Phihh_2(Q) + \dots \Phihh_r(Q) ) = 0 \eeq
or
\beq \label{5.m2} d_H \Phihh_1(Q	)  = 0 , \eeq
\beq \label{5.m3} \delta \Phihh_{j - 1} (Q)
= (-1)^{j+1}  d_H \Phihh_j(Q	)   \in
\Omega^{2r-j+1}_H (G^j), \onebl (1  < j \le r) \eeq
\beq \delta \Phihh_r(Q) = 0. \eeq

\noindent{\em Example.} Let $Q_2 $ be as in   (\ref{4.001}). The
moment $\mu(\theta)$ (corresponding to
the Maurer-Cartan form $\theta$
$ \in \Omega^1(G) \otimes \lieg$,
 which yields the connection $\theta(t)$
on $\bng$ given by (\ref{3.1}))
  is defined by $\mu(\theta)(\phi)  = - \iota_{\tilde{\phi} }
\theta$. Now for the   action of $H = G $ on $G $  by left
multiplication, we have
\beq \label{5.9a} \mu(\theta) (\phi)_g = -  {\rm Ad} (g^{-1} ) \phi \eeq
 (in the notation of
(\ref{4.002})). This leads to
\beq \Phihh_1(Q_2 ) = \Phi_1(Q_2) - \Theta  = - \l - \Theta, \eeq
where \beq \Theta (\phi) = \langle \phi, \theta + \bar{\theta} \rangle
\in \Omega^3_H(G).  \eeq
We also have
\beq \Phihh_2(Q_2 ) = \Phi_2(Q_2) = \Omega =
\langle \theta_1, \bar{\theta}_2 \rangle \eeq
(cf.  (\ref{4.004})).
This result was obtained by Weinstein ([W], Section 4).
(Differences between our formulas and Theorem 4.5 of [W]   result
from Weinstein's use of the right adjoint action on $NG$ rather
than the left adjoint action: see the first sentence of the proof
of [W], Lemma 4.1.)

\renorm
\section{The classifying map for the universal bundle} \label{s6}
In this section, we describe our construction (from Section 2) of the
universal bundle $\pun \to B = \Sigma \times \yyb$ in terms of simplicial
manifolds.

We now take the group $G$ in (\ref{2.1}-\ref{2.2}) to be
the fundamental group $\fund$ of the closed 2-manifold $\Sigma$.
Recall that $K$ was a compact connected semisimple
Lie group and $\hatk = K/Z(K)$.
Then we define simplicial manifolds $\ee$, $\bb$ by
\beq \label{5.1} \ee(n) = (\fund^{n+1} \times \yyb \times
\hatk  )/\fund   \eeq
\beq \label{5.2} \bb(n) = \fund^n \times \yyb,
    \eeq
where $\Pi$ acts by (\ref{1.2}).
Again the face maps of $\ee$ are given by
(\ref{3.6a}), while those
of $\bb$ are given
by (\ref{3.6b}).
There is a (simplicial)
projection map
$\pi_n: \ee(n) \to \bb(n) $ given by
\beq \label{5.3} \pi_n (p_0, \dots, p_n, \rho, k)
= (p_0 p_1^{-1}, \dots, p_{n-1} p_n^{-1}, \rho). \eeq
These maps fit together to give a map of simplicial manifolds
$\ee \stackrel{\pi}{\lrar} \bb$.
The spaces $\ee$ and $\bb$ are equipped with actions of $K$ given
by
\beq \label{5.5} \ttt \in K:
(p_0, \dots, p_n, \rho, k) \mapsto
(p_0, \dots, p_n, \ttt \rho \ttt^{-1}, \ttt k)  \eeq
\beq \label{5.6} \ttt \in K:
(p_1, \dots, p_n, \rho) \mapsto
(p_1, \dots, p_n, \ttt \rho \ttt^{-1})  \eeq
Since $\tsig = E \fund $ and $\Sigma = B \fund$, the simplicial
manifolds $\ee$ and $\bb$ (with $K$ actions) together with the
$K$-equivariant map
$\pi: \ee \to \bb$ are the simplicial realization of the $K$-equivariant
bundle
$\pun \to B = \Sigma \times \yyb$.

\nc{\ev}[1]{{ \rm ev}_{ {#1} } }

We can also define maps of simplicial manifolds
$\psie: \ee \to \bnk$, $ \psib: \bb \to \nk$ for which the diagram
\beq \label{5.4}
\begin{array}{lcr}
\ee & \stackrel{\psie}{\lrar} & \bnk \\
\downarrow & \phantom{\stackrel{\psie}{\lrar} } & \downarrow \\
\bb & \stackrel{\psib}{\lrar} & \nk \\ \end{array}
\eeq
commutes. We take
$$ \psie (n): \ee(n) \to \bnk (n) $$ to be given by
\beq \label{5.05} \psie(n) (p_0, \dots, p_n, \rho, k)
 = ( \ev{p_0} (\rho)  k, \dots, \ev{p_n} (\rho) k) \eeq
while $\psib(n): \bb(n) \to \nk(n) $ is defined by
\beq \label{5.06} \psib(n) (p_1, \dots, p_n, \rho)
 = ( \ev{p_1} (\rho) , \dots, \ev{p_n}(\rho) ) \eeq
Here,
$\ev{p} : {\rm Hom} (\FF, K) \to K$ is the evaluation map
associated to an element $p \in \FF$,  and
we have observed that $\yyb$ is contained in $  {\rm Hom} (\FF, K)$.
If $p_j = x_{r(j)} $ in terms of the standard generators of $\free$, and
we identify $\yyb$ with a subspace of $K^{2g}$ so that a point
$\rho \in \yyb$ becomes identified with
${\bf h } = (h_1, \dots, h_{2g} )$ where $h_j = \rho(x_j)$,
then in this notation we have
\beq \psie(n) (x_{r(0)}, \dots, x_{r(n)}, {\bf h}, k) =
(h_{r(0)}k, \dots, h_{r(n)}k) \eeq
and similarly
\beq  \label{6.*1} \psib(n) (x_{r(1)}, \dots, x_{r(n)}, {\bf h}) =
(h_{r(1)}, \dots, h_{r(n)}). \eeq

Further, the following is easily verified:
\begin{prop} The map $\psie$ is equivariant with respect to the
actions of $K$, where $\ttt \in K$ acts on $\ee$ by (\ref{5.5}),
and on $\bnk$ by left multiplication.
The map $\psib$ is $K$-equivariant, where $\ttt \in K$ acts on $\bb$ by
(\ref{5.6}), and on $\nk$ by the adjoint action.
\end{prop}
Thus we see that $\psib$  is the simplicial
realization of the classifying map for the
bundle $\pun \to \Sigma \times \yyb$ (defined
in Section 2), {\em as a $K$-equivariant bundle}.

\noindent{\em Remark:} Since $\Sigma  = B \fund$ has cohomology
only in dimensions $0,1,2$, we shall in fact only
be concerned with $\nk(n) $ for $n \le 2$.

\renorm
\section{Pullbacks of equivariant characteristic classes} \label{s7}

We now use the classifying maps $\psie$ and $\psib$ from Section 6 to
pull back the equivariant characteristic classes constructed in
Section 5. In this section
we restrict to $K=SU(N)$, so $  \hatk
= PSU(N)$, and  we assume that
$\beta$ is a generator of $Z(K)$. In this situation
$\calm_\beta$ is a smooth manifold.
 This choice of $K$ and $\beta$ is the case for which the generators of
$H^*(\calm_\b)$ are given explicitly in Section 9 of [AB].

Let $Q_r $ $(r \ge 2)$ be the invariant polynomial on
$ \liek = su(N)$ which corresponds to the $r$'th Chern class. Then
according to [AB] (Sections 2 and 9), the generators of  the ring
$H^*(\calm_\beta) $ (using rational, real or complex coefficients) are
given as follows. We fix generators
$c \in H_2 (\Sigma; \ZZ)$ and $\a_j \in H_1(\Sigma; \ZZ) $
$(j = 1, \dots, 2g)$. Let us define a rank $N$ vector bundle
$\uu$ over $\Sigma \times \yyb$ (equipped with a $K$ action covering
the $K$ action on $\yyb$), whose transition functions are
in $SU(N)$ and are lifts of the transition functions of the
principal $\hatk$-bundle $\pun$ in (\ref{2.02}) (or
$\ee$ in  (\ref{5.1})) from $
\hatk = PSU(N)$  to $SU(N)$. (Such a lift exists according to
Section 9 of [AB]; moreover,
the characteristic classes of the vector  bundle $\uu$
defined using such transition
functions  are independent of the choice of lift.)
We then have the
following description of a set of generators of $H^*(\calm_\b)$
$\cong \hk(\yyb)$:
\beq f_r = (c, \crk(\uu) )  \eeq
\beq b_r^j= (\a_j, \crk(\uu) ) \eeq
\beq a_r= (1, \crk(\uu) ) \eeq
Here,  $\crk$ is the $r$'th  equivariant
Chern class, and  $(\cdot, \cdot ) $ denotes the canonical
pairing $H_* (\Sigma ) \otimes \hk(\Sigma \times \yyb) \to
\hk(\yyb) \cong H^*(\mb). $
In the simplicial manifold description of
the universal bundle (\ref{5.1}-\ref{5.2}),
 we may represent the generators
of $H_*(\Sigma)$ by $c \in H_2 (\fund)$, $\a_j \in H_1 (\fund)$
where $H_j (\fund)$ denotes Eilenberg-Mac Lane group homology (see e.g.
[Mac]).
Thus we find the following
\begin{theorem} \label{t7.1} The generators $a_r$, $b_r^j$ and $f_r$
of $\hk(\yyb)$ are given in terms of the map $\Phih$ from Section 5.3
 by the following elements $\tar, \tbrj$ and $\tfr$ of
$\Omega^*_K(\yyb)$:
\beq \tbrj =  \Bigl (\a_j, \psib^* \Phih_1 (Q_r) \Bigr ) \eeq
\beq \tfr = \Bigl ( c, \psib^* \Phih_2 (Q_r) \Bigr ), \eeq
while $\tar$ is the image of  $Q_r \in S(\lieks)^K   $
in $\Omega^*_K(\yyb)$. Here, the classifying map
$\psib$ was defined by (\ref{5.06}).
\end{theorem}

Let us now identify $\yyb$  with a subspace of
$K^{2g}$, and write the projection
on the $j$'th copy of $K$
as $\pi_j: K^{2g} \to K$. We then have the following result.
\begin{prop} \label{p7.*2}
The classes  $f_r$ and
$b_r^j$  are represented in the Cartan model by  $K$-equivariant
polynomial maps
$$ \phi \in \liek \mapsto \tilde{f_r} (\phi), \tilde{b_r^j} (\phi)
 \in \Omega^*(\yyb). $$
Through our choice of generators $x_1, \dots, x_{2g} $ for
$\FF$, we have identified $\yyb$ with a subvariety
of $K^{2g}$.
Let $\theta$ $\in \Omega^1(K) \otimes \liek$
denote  the (left invariant) Maurer-Cartan form.

In terms of this notation, each of the forms
 $\tfr(\phi),
\tbrj(\phi)  \in \Omega^*(\yyb)$ is
 given by a linear combination of polynomials in the quantities
$\ev{p}^* (\theta)   $  and  ${\rm Ad}(\ev{p})(\phi)$.
Here, $\ev{p}: {\rm Hom} (\FF, K) \to K$
is the evaluation map associated to an element $p \in \FF$, and
$p$ ranges over a finite set of elements in $\FF$, while
$j = 1, \dots, 2g$.
\end{prop}
\Proof
For any $ \phi \in \liek$, we have
\beq \label{7.3a} \BPhih{2}(Q_r )(\phi)  =
\sum_{A, B} \k_{AB} Q_r\Bigl (\theta_{\a_1}, \theta_{\a_2},
[\theta_{\a_3},\theta_{\a_4} ], \dots,
[\theta_{\a_{2r-2m- 3} },\theta_{\a_{2r-2m-2} } ],
{\rm Ad} (k_{\b_1}^{-1} ) \phi, \dots, {\rm Ad} (k_{\b_m}^{-1} ) \phi
\Bigr )  \eeq
$$ \phantom{aaaaaaa}\in \Omega^{2r-2-2m} (K^3). $$
Here, $A = (\a_1, \dots, \a_{2r-2-2m})$,
$B = (\b_1, \dots, \b_m) $ are multi-indices, and
$K_{AB} $ are certain constants.  Also, $\a_i, \b_i $ $ = 0, 1, 2,$
and $k_\b$ refers to the $\b$-th copy of $K$ in $K^3$.
(Here, the factors involving $\phi$  result from  (\ref{5.9a}).)

Thus also by  (\ref{4.2a}) we have
\beq \label{7.4c} \Phih_2(Q_r) (\phi) =
\sigma_2^* \BPhih{2}(Q_r) (\phi) \in \Omega^{2r-2-2m} (K^2), \eeq
where
$\sigma_2: K^2 \to K^3$ is the map
\beq \label{7.4b}
\sigma_2: (k_1, k_2) \mapsto (1, k_1, k_1 k_2)  \eeq

So we have
\beq \label{7.5a}
\tfr(\phi) = \Bigl (c, \psib^* \Phih_2(Q_r) (\phi) \Bigr )
 \in \Omega^{2r-2-2m}
\Bigl ({\rm Hom} (\FF, K) \Bigr ) \eeq
Finally we use the explicit description of $c \in C_2 (\free)$
([G], above Proposition 3.1):
\beq \label{7.*1}
c = \sum_{j = 1}^{2g} (\frac{\partial R}{\partial x_j} \otimes x_j), \eeq
where $\partial R/\partial x_j $ $ \in C_1(\free)$ denotes the
differential in the Fox free differential calculus (see e.g. [G], Sections
 3.1-3.3).
We may write
$\partial R/\partial x_j = \sum_f a_f f$ where $a_f \in \ZZ$ and
$f$ are certain elements in $\free$ given explicitly by words in
the generators $x_1, \dots, x_{2g}$. We have in fact
\beq \frac{\partial R}{\partial x_j} =
\g_j^0 - \g_j^1, \eeq
where $  \g_j^\a \in \FF$   is given by
\beq \begin{array}{lcr}
\g_{2j-1}^0 &=& \prod_{l = 1}^{j-1} [x_{2l-1}, x_{2l} ]  \\
\g_{2j-1}^1 &=& \g_{2j-1}^0 x_{2j-1} x_{2j} x_{2j-1}^{-1} \\
\g_{2j}^0 &=& \g_{2j-1}^0 x_{2j-1} \\
\g_{2j}^1 &=& \g_{2j-1}^0 [x_{2j-1}, x_{2j} ] . \end{array} \eeq
Thus we have from (\ref{7.5a})
\beq \label{7.5b} \tfr(\phi) =
\sum_{j = 1}^{2g} \sum_{\tau = 0 }^1  (-1)^\tau
\Bigl (\ev{\g_j^\tau} \times \ev{x_j}  \Bigr )^* \Phih_2 (Q_r) (\phi) \eeq
where
$\ev{p} : {\rm Hom} (\FF, K) \to K $ is the evaluation map associated
to an element $p$ in $\FF$ and one defines
$\ev{\g_j^\tau} \times \ev{x_j}: {\rm Hom} (\FF, K) \to K^2$.

Combining (\ref{7.5b})  with  (\ref{7.4c}) we see  that
\beq \label{7.5c}
\tfr(\phi)
= \sum_{j= 1}^{2g} \sum_{\tau = 0}^1  (-1)^\tau
 \Psi_{j,\tau}^* \BPhih{2}(Q_r)(\phi) , \eeq
where
$\Psi_{j,\tau} = \ev{1} \times \ev{\g_j^\tau} \times \ev{\g_j^\tau x_j} : $
$  {\rm Hom}(\FF, K) \to K^3$ is the map
\beq \Psi_{j,\tau}: \rho \mapsto (1, \rho(\g_j^\tau), \rho(\g_j^\tau x_j)).
\eeq
Using the formula (\ref{7.3a}) for
$ \BPhih{2}(Q_r) (\phi)$, this gives the desired result.

Explicitly,
if we decompose the right
hand side of (\ref{7.3a}) as
$\BPhih{2}(Q_r)(\phi)  = \sum_{a}  \kappa_a
 \omega^0_{a}\wedge  \omega^1_{a}\wedge $
$ \omega^2_{a} $
for
$\omega^j_{a}$ $\in \Omega^*(K_j) $ (where
$K_j$ is the $j$-th copy of $K$ in $K^3$), we have
\beq \label{7.6a}
(c, \psib^* \Phih_2(Q)(\phi)  )
= \sum_{j= 1}^{2g} \sum_{\tau = 0}^1  (-1)^\tau
  \sum_{a}  \kappa_a (\ev{1}^* \omega^0_a)
(\ev{\g^\tau_j}^* \omega^1_a) (\ev{\g^\tau_j x_j}^* \omega^2_a). \eeq
Let us define elements $z_{j,\a}^\tau$ in $\FF $ (for
$ \a = 0,1,2$ and $\tau = 0, 1$) by
$z_{j,0}^\tau = 1, $ $ z_{j,1}^\tau = \g^\tau_j, $
$z_{j,2}^\tau = \g^\tau_j x_j$.
Then we have using (\ref{7.3a})  and (\ref{7.5a})
\beq \label{7.6b}
\tfr (\phi) =
\sum_{j = 1}^{2g} \sum_{\tau = 0}^1 (-1)^\tau
\sum_{A, B} \k_{AB} Q_r \Bigl (
\ev{z_{j,\a_1}^\tau }^* \theta, \ev{z_{j,\a_2}^\tau }^* \theta,
[\ev{z_{j,\a_3}^\tau}^* \theta, \ev{z_{j,\a_4}^\tau}^* \theta ], \dots,
\eeq
$$ \phantom{aaaaaaa}\dots [ \ev{z_{j,\a_{2r-2m-3}}^\tau }^* \theta,
\ev{z_{j,\a_{2r-2m-2}}^\tau}^* \theta ],
{\rm Ad} (\ev{(z_{j,\b_1}^\tau)^{-1} })( \phi  ) , \dots,
{\rm Ad} (\ev{(z_{j,\b_m}^\tau)^{-1} })( \phi  ) \Bigr )
\in \Omega^{2r-2-2m}(\yyb). $$

 The proof for $\tbrj(\phi) $ is
similar but easier. $\square$

\noindent{\em Remark:} If $m: K \times K \to K$ is the
 multiplication map,
and
  $i: K \to K$ is the inversion map, we have
 \beq
 m^* \theta  = \theta_2 + {\rm Ad}(k_2^{-1}) \theta_1, \bla
(i^* \theta)_k  = {\rm Ad}(k^{-1}) \theta. \eeq
Applying this to    (\ref{7.6b})  we get an alternative
expression
in terms of the generators $x_j$ for $\FF$.

\renorm
\section{Extended moduli spaces} \label{s8}

The
 {\em extended moduli space } $\xb$   was defined in
[J1];
 it   is a symplectic space equipped with a Hamiltonian action of
$K$,
such that the
space $\yyb$ treated above embeds in
$\xb$ as the zero locus of the moment map.
Thus the symplectic reduction of $\xb$
is  the space $\mb = \yyb/K$ associated
to the central element $\b$ $\in Z(K)$. The symplectic form on $\xb$ was given
a gauge
theoretic construction in [J1], and a construction via group
cohomology (using techniques similar to the present paper) was
given in [J2] and independently in [H].
 In this section we shall extend the results of the previous
sections to obtain de Rham representatives
for classes in $\hk(\xb) $ whose restrictions to
$\hk(\yyb)$ are the classes constructed in Section 6.

Let us first recall from [J1] or [J2] the definition of the
space $\xb$. We had defined an element $
R = \prod_{j = 1}^g x_{2j-1} x_{2j} x_{2j-1}^{-1}
x_{2j}^{-1} \in \free
 $ in terms of the generators $x_j$ of   $\free$.
The space $\xb$ is defined
as a fibre product
\beq   \label{8.1}
\xb = (\epsr \times \epmb)^{-1} (\diag) \subset K^{2g} \times \liek. \eeq
Here,
${\rm Hom} (\FF, K)$ has been identified with $K^{2g}$, while
$\diag$ is the diagonal in $K \times K$ and
$\epsr: K^{2g}  \to K$ is the map $\epsr: \rho \in
{\rm Hom} (\FF, K) \mapsto \rho(R) . $ The map
  $\epm: \liek \to K$ is the
exponential map, and $\epmb = \beta \cdot \epm$ in terms
of the element $\beta \in Z(K)$.
The space $\xb$ is equipped with two canonical
projection
maps $\proj_1: \xb \to K^{2g}$ and $\proj_2: \xb \to \liek$, and there is
a commutative diagram
\beq \label{7.2}
\begin{array}{lcr}
\xb  & \stackrel{\proj_2}{\lrar} & \liek \\
\scriptsize{\proj_1}
\downarrow & \phantom{\stackrel{\psie}{\lrar} } & \downarrow \scriptsize{e_\b}
 \\
K^{2g}  & \stackrel{\epsr }{\lrar} & K\\ \end{array} \eeq

Define  $\sio = \Sigma - D^2$. We shall now adapt the techniques from
Section \ref{s6} to construct a principal $ \hatk$ bundle $\eex$
over $\sio \times \xb$, equipped with a distinguished homotopy
class of trivializations on $\partial \sio \times \xb$, and with an
action of $K$ covering the action of $K$ on $\xb$. This $\hatk$ bundle
will restrict on $\yyb$ to the $K$-equivariant
bundle $\ee$ constructed in Section \ref{s6}; hence its
equivariant characteristic classes are the extensions to $\xb$ of the
equivariant characteristic classes of the bundle $\ee$, which were used in
Section \ref{s7} to construct the generators of the cohomology
ring $H^*(\mb)$.

We start with a construction of spaces
analogous to (\ref{5.1}-\ref{5.2}) and maps analogous to
(\ref{5.05}-\ref{5.06});  the spaces
(obtained from an action of $\free$ similar to (\ref{1.2})) are
 \beq \label{8.01f} \eex(n) \eqdef
 ( \free^{n+1} \times \xb \times \hatk)/\free , \eeq
\beq \label{8.02f}   \bbx(n) \eqdef \free^n \times \xb.
  \eeq
These bundles are pulled back in an obvious way
(via the maps $\proj_1: \xb \to K^{2g} $) from bundles
 \beq \label{8.01fn} \ceex(n) \eqdef
 ( \free^{n+1} \times K^{2g} \times \hatk)/\free , \eeq
\beq \label{8.02fn}   \cbbx(n) \eqdef \free^n \times K^{2g}.
  \eeq
Denote the canonical maps $\eex \to \ceex$
(resp.  $\bbx \to \cbbx$)  by $p_E$
(resp. $p_B$): then there is an obvious commutative diagram
\beq \label{7.2*1}
\begin{array}{lcr}
\eex  & \stackrel{p_E}{\lrar} & \ceex \\
\phantom{\scriptsize{\proj_1}}
\downarrow & \phantom{\stackrel{\psie}{\lrar} } & \downarrow
\phantom{\scriptsize{e_\b} }
 \\
\bbx  & \stackrel{{p_B}}{\lrar} & \cbbx \\ \end{array} \eeq

The maps
 \beq \label{8.01}\psiex(n): \eex(n) =
 ( \free^{n+1} \times \xb \times \hatk)/\free  \to \hatk^{n+1}, \eeq
\beq \label{8.02} \psibx(n) :  \bbx(n)  = \free^n \times \xb
  \to \hatk^n \eeq
are  given by $\psiex =  \cpsiex \circ p_E$, $\psibx =
\cpsibx \circ p_B $, where
$$\cpsiex(n):
( \free^{n+1} \times  {\rm Hom}(\free, K) \times \hatk)/\free  \to \hatk^{n+1},
$$
$$ \cpsibx(n): \free^n \times  {\rm Hom}(\free, K)
 \to \hatk^n$$ are
defined (as in (\ref{5.05}-\ref{5.06})) by
\beq \label{8.05p} \cpsiex(n) (p_0, \dots, p_n, \rho, k)
 = ( \rho(p_0) k, \dots, \rho(p_n) k) \eeq
\beq \label{8.06p} \cpsibx(n) (p_1, \dots, p_n, \rho)
 = ( \rho(p_1) , \dots, \rho(p_n)). \eeq
 These maps give a classifying
map for the bundle
\beq  \eex = (E \free \times \xb\times  \hatk  )/\free
 \to \bbx =  B \free \times \xb. \eeq
In other words, the following diagram
commutes:
\beq  \label{7.1}
\begin{array}{lcr}
\eex  & \stackrel{\psiex}{\lrar} & \bnk \\
\downarrow & \phantom{\stackrel{\psie}{\lrar} } & \downarrow \\
\bbx & \stackrel{\psibx}{\lrar} & \nk \\ \end{array} \eeq

Of course, $\sio$ has the homotopy type of $B \free$, so
$\eex$ may be regarded as a bundle over $\sio \times \xb$. The action of
$K$ on $\eex$ is given by the following action on
$E \free \times \xb \times \hatk$:
\beq \label{8.z2} \ttt \in K: \; (e, (\rho,\L), k)\mapsto
(e, (\ttt \rho \ttt^{-1},\ttt \L \ttt^{-1}), \ttt k). \eeq
(Here, $(\rho, \L) \in \xb \subset K^{2g} \times \liek ;$  see (\ref{8.1}).)
The maps in the  diagram (\ref{7.1})
are $K$-equivariant, where $K$ acts
on $\bnk$ by left multiplication
and on $\nk$ by the adjoint action.

We now construct a trivialization of $\eex$ over $\partial \sio
\times \xb$.
Consider the map $\psibx(1): \onebl \bbx(1) = $
 $\free \times \xb
 \to  \nk(1) = \hatk $ determined by (\ref{7.1}).
Its restriction
to $\ZZ \cdot R \times \xb$ (where $\ZZ \cdot R$ is the
abelian subgroup of $\free$ generated by the element $R$)
lifts to a map
$\ZZ \cdot R \times \xb \to \liek$, because of the diagram (\ref{7.2}).
Since $\liek$ is contractible, this implies there is a canonical
trivialization of the bundle $\eex$ over $\ZZ \cdot R \times \xb$.
Now since $\ZZ \cdot R = \ZZ = \pi_1 (\partial \sio) $, and
in fact the map
$\ZZ \cdot R \to \free$  gives rise to a map $B \ZZ  \to B \free$
which is identified with the map $\partial \sio \to \sio$, we have
constructed a trivialization of the bundle $\eex$ over
$\partial \sio \times \xb$.

Via  the classifying map given by (\ref{7.1}),
the canonical $K$-invariant connection $\theta(t) $
(\ref{3.1}) on  $\bnk \to \nk$
 pulls back to a $K$-invariant connection on  the bundle
$\eex \to \sio \times \xb$. The equivariant characteristic
classes of $\eex$ are then given by pulling back the classes\footnote{Recall
that $ K$ acts on $\bnk(n) = \hatk^{n+1} $ by left multiplication.}
$\Phih(Q_r) $
$\in H^*_K(\nk) $
 via the classifying map $\psibx$. Since the restriction
of $\psibx$ to $\partial \sio \times \xb$ factors through
a map to a contractible space, the classes
$\psibx^* \Phih ( Q_r ) $ live in $H^* (\sio, \partial \sio) \otimes
\hk(\xb)$. Their restrictions to $(\sio, \partial \sio)
\times \yyb$ identify
(since $H^*(\sio, \partial \sio) $ $\cong H^*(\Sigma)$)
with the classes
in $H^*(\Sigma) \otimes \hk(\yyb)$ constructed in Section \ref{s7}.

The generators $\tar$, $\tbrj$
and $\tfr$ of $H^*(\mb) \cong
\hk(\yyb)$ were obtained in Section \ref{s7} by
 pairing ${\psib}^* \Phih(Q_r) $ with the generators of $H_*(\Sigma)$
via the slant product.
We obtain the equivariant extensions of $\tar, \tbrj $ and  $ \tfr$ to
$\xb$ by pairing
$\psibx^* \Phih(Q_r) $ with the generators of
$H_*(\sio, \partial \sio) \cong H_*(\Sigma)$.
To construct these equivariant extensions explicitly,  we must
introduce the equivariant complex $\Omega^*_K(\liek)$, where $ K$
acts on $\liek$ via the adjoint action. Now there is an equivariant
map $\hop: \Omega^{*+1}_K
(\liek) \to \Omega^*_K (\liek)$ which appears in
the usual proof of the
Poincar\'e lemma (see e.g.  [Wa],  Lemma 4.18). The
definition is as follows:  for any $v \in
\liek$ and $\beta \in \Omega^*_K(\liek)$ we define
\beq
\label{p.2} (\hop \beta)_v =
\int_0^1 F_t^* (\iota_{\bar{v} } \beta ) dt \eeq
where $F_t: \liek \to \liek$ is the map given by multiplication
by $t $ and
 $\bar{v} $ is the vector field on $\liek $ which takes the constant
value $v$.
The map $\hop$ is usually defined as a chain homotopy in the
ordinary de   Rham complex on a vector space; it
generalizes to the equivariant de Rham complex on the
vector space $\liek$ because the
maps $F_t$ are $K$-equivariant maps.  More precisely, the map $\hop$
 satisfies the  property
\beq \label{8.001}  \hop d_K + d_K\hop = 1,  \eeq
in other words $\hop$ defines a  chain homotopy equivalence
 (which is in fact
equivariant).
The formula
(\ref{8.001})
  follows because  $\hop d + d \hop = 1$  and also
$\hop \iota_\phit = - \iota_\phit \hop$, where
$\iota_\phit $ is the interior product with the
vector field $\phit$ generated by the action of an
element $\phi \in \liek$.

Let us denote
the generator of $H_2(\sio, \partial \sio) $ $\cong H_2(\free, \ZZ
\cdot R) $ by $c$ and those of $H_1(\sio, \partial \sio) $
$ \cong H_1(\free, \ZZ \cdot R) $ by $\a_j$
($j = 1, \dots, 2g$). We then have

\begin{theorem}  \label{t8.1} Let $K = SU(n)$.
The extensions to $\hk(\xb)$ of the generators
 $\tbrj, \tfr$  of
$  \hk(\yyb)$  are  represented by equivariant differential forms
$\tar(X), \tbrj(X), \tfr(X)  \in  \Omega^*_K(\xb) $ defined by
\beq \label{8.f}
 \tfr(X) = \proj_1^*\tfr(X)_1   + \proj_2^*
\tfr(X)_2  \eeq
\beq \label{8.g}
\tbrj(X) = \proj_1^* \tbrj(X)_1 \eeq
where
\beq \label{8.17a}
 \tfr(X)_1 = (c, \cpsibx^* \Phih_2 (Q_r) ),    \eeq
\beq \label{8.17b} \tfr(X)_2 = - \hop \Bigl (
e_\b^* \Phih_1 (Q_r) \Bigr ),  \bla \mbox{and} \eeq
\beq  \label{8.b}
 \tbrj(X)_1 = (\a_j, \cpsibx^* \Phih_1(Q_r) ). \eeq
The equivariant differential form  $\tar(X)$  representing the extension
of $\tar$ to $\Omega^*_K(\xb)$ is   the image of
$Q_r \in S(\lieks)^K  $ in $\Omega^*_K(\xb)$.
Here, $ K $ acts on $\liek$ by the adjoint action.
Also,  $\hop: \Omega^{p+1}_K
(\liek) \to \Omega^p_K (\liek)$ is the
chain homotopy map  defined in (\ref{p.2}), and
$Q_r$ is the  invariant polynomial on $\liek$ corresponding to the
$r$th Chern class.
\end{theorem}

\Proof
 Let
$Q$ be an invariant polynomial of degree
$r$ on $\liek$.
The class $\Phih (Q ) $ $ = \{ \Phih_n(Q) \} $
$ \in \Omega^*_K(\nk) $ satisfies
$ \Phih_1 (Q) \in \Omega^{2r-1}_K(\hatk)$,
$ \Phih_2 (Q) \in \Omega^{2r-2}_{\hmac}(\hatk^2)$. Since
$\Phih (Q) $ is closed under the total differential
$d_{\hmac}  \pm \delta$\footnote{Recall the sign convention for
total differentials of double complexes introduced in Section 4.}
   of the complex
$\cala^{*,*}_{\hmac} (\nk)$, we have
\beq \delta
\Phih_1(Q) =   -
d_{\hmac} \Phih_2 (Q) \in  \Omega^{2r-1}_{\hmac}(\hatk^2). \eeq
Pulling back via the classifying map $\psibx$, we find
\beq  d_{\hmac} \Bigl ( \psibx^* \Phih_2 (Q) \Bigr )
= \psibx^*(d_{\hmac} \Phih_2 (Q) )  \eeq
$$ \bla =
- \psibx^* \delta  \Phih_1 (Q) =
 -  \delta \psibx^* \Phih_1 (Q). $$
We now consider the equivariant complex $\Omega^*_{\hmac}(\liek)$, where
$ K$ acts on $\liek$ by the adjoint action.
Now  $e_\b^* \Phih_1(Q) \in \Omega^{2r-1}_{\hmac} (\liek) $ is equivariantly
closed (by (\ref{5.m2})), so
since $H^*_{\hmac}(\liek) $ $\cong H^*(BK)$ is nonzero only
in even dimensions,   there is
an element
$\sigma_Q \in \Omega_{\hmac}^{2r-2} (\liek) $ such that
$e_\b^* \Phih_1 (Q) = d_{\hmac} \sigma_Q$.
 In fact, because of (\ref{8.001})  one may choose
\beq  \label{8.z3} \sigma_Q = \hop (e_\b^* \Phih_1(Q)). \eeq
The idea
of using  the operator $\hop$ to construct closed forms
is due originally
to Weinstein, who suggested this approach in the case $Q = Q_2$.

We thus find that
$d_{\hmac} (\psibx^* \Phih_2 (Q) ) =  -\delta \psibx^* \Phih_1(Q)$.
So if $c \in C_2 (\free)$ is the element representing the
generator of $H_2(\sio, \partial \sio)$, for which
$\partial c = 1-R \in C_1(\free)$ ([G], above Proposition 3.9),
then we have
\beq d_{\hmac} \Bigl (c, \psibx^* \Phih_2 (Q) \Bigr )  =
-\Bigl (c, \delta \psibx^* \Phih_1 (Q) \Bigr )
\eeq
$$ \onebl = - (\partial c, \psibx^* \Phih_1 (Q) ) $$
$$ \onebl = -(1-R, \psibx^* \Phih_1 (Q)) $$
$$ \onebl =  (\epsr \circ \proj_1)^* \Phih_1 (Q) \in \Omega^*_K(\xb)
\onebl (\mbox{cf. (\ref{7.5b})}) $$
$$ \onebl =   \proj_2^* (e_\b^* \Phih_1(	Q) )
\onebl \mbox{by (\ref{7.2})}. $$
Thus because of (\ref{8.001}) and (\ref{8.z3})  we easily obtain
$$d_{\hmac}(c, \psibx^* \Phih_2 (Q) ) = d_{\hmac} \proj_2^*( \sigma_Q), $$
so that  the element
$$ \proj_1^*  (c, \cpsibx^* \Phih_2 (Q) ) - \proj_2^*( \sigma_Q)  $$
is equivariantly closed in $\Omega^{2r-2}_{\hmac}(\xb)$.
When  $Q = Q_r$ is
the invariant polynomial on $\liek = su(N)$
corresponding to the $r$'th Chern class,
this element is  the extension to $ \Omega^*_K(\xb) $ of the
element $\tfr$  described in Section 7.

If $\a_j$ ($j =1, \dots, 2g$) are the generators of
$H_1 (\sio) = H_1 (\free)$,  we find that
$$d_{\hmac} (\a_j, \psibx^* \Phih_1 (Q) ) = 0 $$
since $d_{\hmac} \Phih_1(Q) = 0 $. The element
$(\a_j, \psibx^* \Phih_1 (Q) )$ represents the extension
of the element $\tbrj \in \Omega^{2r-1}_K(\yyb)$ to
$\Omega^{2r-1}_K(\xb)$. $\square $

\noindent{\em Remark.} The appearance of a second term in (\ref{8.f})
corresponding to the boundary of $\sio$ should not be surprising.
We know a  de Rham representative of a differential form $\psibx^*
\Phih_2(Q)$ lives in $\Omega^*(\sio) \otimes \Omega^*_K(\xb)$ and
is equivariantly closed. Evaluating such a form on $c \in C_2(\free)$
corresponds to integrating it over $\sio$. Now for any
Cartesian product of  manifolds
$S \times Z$ with $\partial Z = 0 $, Stokes' Theorem for integration
over the fibre $S$ takes the form
\beq \int_S d_{S \times Z} \alpha - (-1)^{{\rm dim} S} d_Z \int_S \alpha
+ \int_{\partial S} \alpha = 0. \eeq
The same is true when $Z$ is a closed manifold acted on by a group
$K$ (and $K$ acts trivially on $S$), when the de Rham exterior
differential $d$ is replaced by the
Cartan model exterior differential $d_K$.
Taking $Z = \xb$ and  $S = \sio$, we obtain when $\alpha $ is the
equivariantly closed differential form
$\psibx^*
\Phih_2(Q) $,
\beq \label{8.0z1} d_K \int_\sio \psibx^*
\Phih_2(Q) = \int_{\partial \sio} \psibx^*
\Phih_2(Q). \eeq
But over $\partial \sio \times \xb$, the classifying map
$\psibx$ factors through an  equivariant map to a contractible
space (see the paragraph after (\ref{8.z2})), so the right hand side of
(\ref{8.0z1}) can explicitly be written as a $d_K$-exact form. One choice
of such an exact form is the image under $d_K$ of
 the second term on the right hand side of  (\ref{8.f}).

\noindent{\em Example.} Let $K = SU(n)$, and let
$Q_2$ $ = \langle \cdot , \cdot
\rangle $ be the invariant inner product on $\liek$ that gives
rise to the second Chern class. The corresponding generator
$\tilde{f}_2$
is the cohomology
class of the standard  symplectic
form on $\calm_\beta$. In this case,
 the formula  (\ref{8.f}) generalizes the result of
[J2], which gives the extension of the symplectic form $\omega$ on
$\calm_\beta = \yyb/K$ to a closed 2-form $\tilde{\omega}$
on $\xb$ (indeed a symplectic form
on a neighbourhood of $\yyb$ in $\xb$). This 2-form is defined by ([J2], (5.5))
\beq \tilde{\omega} = \proj_1^* \omega - \proj_2^* \sigma, \eeq
where $\sigma = \sigma_{Q_2}   = \hop (e_\beta^* \Phi_1(Q_2) )$
and $$ \Phi_1(Q_2) = -\lambda =
-(1/6) \langle \theta, [\theta, \theta] \rangle
  \in \Omega^3 (\hatk). $$
The 2-form $\omega \in \Omega^2 (K^{2g})$  is defined
(see (\ref{4.004})) by
\beq \omega = (c, \psibx^* \Omega) \eeq
where $\Omega = \Phi_2 (Q_2)  =
\Phi_2^K(Q_2) = \langle \theta_1, \bar{\theta_2} \rangle
\in \Omega^2 (K^2). $
(Here, $\theta_i$ is the left invariant Maurer-Cartan form on the
$i$'th copy of $K$, while $\bar{\theta_i} $ is the corresponding
right-invariant form.)

The corresponding {\em equivariantly}
closed 2-form $\bar{\omega} $
on $\xb$ is  the element in $\Omega^2_K(\xb) $ given
in terms of the Cartan model  (for $\phi \in \liek$) by
\beq \bar{\omega} (\phi) = \tilde{\omega} + \langle \mu, \phi \rangle, \eeq
where $\mu = - 2 \proj_2: \xb \to \liek$ is  the moment map for
the action of $K$ on $\xb $ ([J2], Proposition 5.4).
Indeed we have as a special case of (\ref{8.f})
\beq \bar{\omega} = \proj_1^* \omega - \proj_2^* \sigma^K \eeq
where $\sigma^K = \hop e_\beta^*(\Phih_1(Q_2) ). $
Here,
$\Phi_1^K(Q_2) = -\lambda - \Theta$,
where $\Theta \in \Omega^3_{\hmac}(K) $ is given by $\Theta(\phi) =
\langle \phi, \theta + \bar{\theta} \rangle $: see the end of Section 5.3.

Finally, we have the following analogue of Proposition \ref{p7.*2}, whose
proof is very similar to the proof of that Proposition.

\noindent{\em Remark:} The equivariant
forms $  \tar(X) $  are (by Theorem \ref{t8.1}) simply the
invariant polynomials $Q_r \in S (\lieks)^K$
$\cong H^*_K ({\rm pt})$; this simple observation
about  the $\tar(X) $ is replaced by the result of
Proposition \ref{p8.2} for the equivariant forms $\tbrj(X)$ and
$\tfr(X)$.

\begin{prop} \label{p8.2}{\bf (a)}
{}~In the notation of Theorem \ref{t8.1}, the
 differential forms    $\tfr(X)_1  $ and
$\tbrj(X)_1  $  are given by
polynomial maps
$$ \phi \in \liek \mapsto \tfr(X)_1  (\phi), \tbrj(X)_1  (\phi)
 \in \Omega^*({\rm Hom} (\FF, K) ). $$
Let  $\theta$ $\in \Omega^1(K) \otimes \liek$
denote  the (left invariant) Maurer-Cartan form; then
each of the forms
 $\tfr(X)_1  (\phi),
\tbrj(X)_1 (\phi)  \in \Omega^*( {\rm Hom}(\FF,K) ) $ is
 given by a linear combination of polynomials in the quantities
$\ev{p}^* (\theta)   $  and  ${\rm Ad}({\ev{p}})(\phi)$.
Here, $\ev{p}: {\rm Hom} (\FF, K) \to K$
is the evaluation map associated to an element $p \in \FF$, and
$p$ ranges over a finite set of elements in $\FF$, while
$j = 1, \dots, 2g$.

\noindent{\bf (b)} The differential forms
$\tfr(X)_2 (\phi) $ from  Theorem \ref{t8.1} are
given in terms of a basis for $\Omega^*(\liek) $
(determined by an orthonormal basis for $\frak{t} $ and a
choice of roots for the action of $K$ on $\frak{t}^\perp$) by
a collection of smooth bounded functions on $\liek$.
\end{prop}

\vspace{0.2in}

{\Large \bf References}

[AB] Atiyah, M.F. and Bott, R., The Yang-Mills equations over Riemann surfaces,
{\em Phil. Trans. R. Soc. Lond.} {\bf A 308} (1982), 523-615.

[BGV] Berline, N., Getzler, E., Vergne, M. {\em Heat Kernels and
Dirac Operators}, Springer-Verlag
(Grundlehren vol. 298), 1992.

[B1]
Bott, R., On the Chern-Weil homomorphism and the continuous cohomology
of Lie groups, {\em Advances in Math.} {\bf 11} (1973), 289-303.

[BSS] Bott, R., Shulman, H., Stasheff, J.,  On the de Rham theory  of
certain classifying spaces, {\em Advances in Math.} {\bf 20 }
(1976) 43-56.

[Du] Dupont, J.L., Simplicial De Rham cohomology and characteristic
classes of flat bundles, {\em Topology} {\bf 15} (1976) 233-245.

[Ge] Getzler, E.,  The equivariant Chern character for non-compact  Lie
groups, {\em Advances in Math.}, to appear.

[G]  Goldman, W., The symplectic nature of fundamental groups of surfaces,
{\em Advances in Math.} {\bf 54} (1984), 200-225.

[H] Huebschmann, J., Symplectic forms and Poisson structures of
certain moduli spaces, preprint (1993).

[J1] Jeffrey, L.C., Extended moduli spaces of flat connections
on Riemann surfaces, {\em Math. Annalen} {\bf 298} (1994)
667-692.

[J2] Jeffrey, L.C., Symplectic  forms on moduli spaces of flat connections
on 2-manifolds,
to appear in Proceedings of the Georgia
International Topology Conference (Athens, GA, 1993), ed. W. Kazez.

[JK] Jeffrey, L.C., Kirwan, F.C., Localization for nonabelian
group actions, preprint alg-geom/9307001 (1993); {\em Topology} (in press).

[JW] Jeffrey, L.C., Weitsman, J. Toric structures on the
moduli space of flat connections on a Riemann surface: volumes
and the moment map, {\em Advances in Math.} {\bf 109}
(1994) 151-168.

[K]
Karshon, Y., An algebraic proof for the symplectic structure of moduli space,
{\em Proc. Amer. Math. Soc.} {\bf 116} (1992) 591-605.

[M] May, J.P., {\em Simplicial Objects in Algebraic Topology},
Van Nostrand (1967).

[Mac] Mac Lane, S., {\em Homology}, Springer-Verlag
 (Grundlehren v. 114) (1963).

[Se1] Segal, G. Classifying spaces and spectral sequences,
{\em Publ. Math. IHES} {\bf 34} (1968) 105-112.

[Se2] Segal, G. Categories and cohomology theories, {\em Topology}
{\bf 13} (1974) 293-312.

[Sh]
Shulman, H.B., Characteristic classes and foliations, Ph.D. Thesis,
University of California, Berkeley (1972).

[T] Thaddeus, M., Conformal field theory and the cohomology of
the moduli space of stable bundles, {\em J. Diff. Geo.}
{\bf 35} (1992) 131-149.

[Wa] Warner, F.W., {\em Foundations of Differentiable Manifolds
and Lie Groups}, Springer-Verlag, 1983.

[W] Weinstein, A., The symplectic structure on moduli space, to appear
in Andreas Floer memorial volume, Birkh\"auser.

\end{document}